\pdfoutput=1
\documentclass[sn-basic]{sn-jnl}
\jyear{2021}%
\usepackage{natbib}
\usepackage{graphicx} 
\usepackage{subfigure}
\theoremstyle{thmstyleone}%
\newtheorem{theorem}{Theorem}%  meant for continuous numbers
\newcommand{\st}{\mbox{subject to}}

\theoremstyle{thmstyletwo}%
\newtheorem{example}{Example}%
\newtheorem{remark}{Remark}%

\theoremstyle{thmstylethree}%
\newtheorem{definition}{Definition}%
\newtheorem{lemma}{Lemma}
\usepackage{natbib}
\begin{document}

\title[A Weighted Randomized Sparse Kaczmarz Method for Solving Linear Systems]{A Weighted Randomized Sparse Kaczmarz Method for Solving Linear Systems}


\author[1]{\fnm{Lu Zhang}}\email{zhanglu21@nudt.edu.cn}

\author[2]{\fnm{Ziyang Yuan}}\email{yuanziyang11@nudt.edu.cn}

\author[1]{\fnm{Hongxia Wang}}\email{wanghongxia@nudt.edu.cn}

\author*[1]{\fnm{Hui Zhang}}\email{h.zhang1984@163.com}

\affil*[1]{\orgdiv{Department of Mathematics}, \orgname{National University of Defense Technology}, \orgaddress{\city{Changsha}, \postcode{Hunan 410073},\country{China}}}

\affil[2]{\orgname{Academy of Military Science}, \orgaddress{\city{Beijing},\country{China}}}


\abstract{The randomized sparse Kaczmarz method, designed for seeking the sparse solutions of the linear systems $Ax=b$, selects the $i$-th projection hyperplane with likelihood proportional to $\|a_{i}\|_2^2$, where $a_{i}^T$ is $i$-th row of $A$. In this work, we propose a weighted randomized sparse Kaczmarz method, which selects the $i$-th projection hyperplane with probability proportional to $\lvert\langle a_{i},x_{k}\rangle-b_{i}\rvert^p$, where $0<p<\infty$, for possible acceleration. It bridges the randomized Kaczmarz and 
greedy Kaczmarz by parameter $p$. Theoretically, we show its linear convergence rate in expectation  with respect to the Bregman distance in the noiseless and noisy cases, which is at least as efficient as the randomized sparse Kaczmarz method. The superiority of the proposed method is demonstrated via a group of numerical experiments.}




\maketitle
\section{Introduction}
The most fundamental problem in linear algebra and numerical mathematics may be to find solutions of the linear systems
\begin{equation}\label{equation:1}
Ax=b,
\end{equation}
with $A\in \mathbb{R}^{m\times n}$ and $b\in \mathbb{R}^m$ being given. Among many existing iterative methods, we focus on Kaczmarz-type methods. Denote the rows of $A$ by $a_1^T, \cdots, a_m^T\in\mathbb{R}^n$ and the entries of $b$ by $b_1,\cdots, b_m\in \mathbb{R}$; then the original Kaczmarz method \citep{1937Angen} selects the $i$-th row cyclically and iterates $x_k$ via
\begin{equation}\label{equation:2}
x_{k+1}=x_k-\frac{\langle a_i,x_k\rangle-b_i}{\|a_i\|_2^2}a_i,
\end{equation}
which means that $x_{k+1}$ is computed as an orthogonal projection of $x_k$ onto the $i$-th hyperplane $H(a_i,b_i):=\{x\in \mathbb{R}^n\lvert~\langle a_i,x\rangle=b_i\}$. The Kaczmarz method becomes more and more popular recently due to the seminal work \citep{2009A}, where the first elegant convergence rate was obtained by considering its randomized variant. Specifically, the randomized Kaczmarz method (RK) selects a projection hyperplane with the probability $p_i=\frac{\|a_i\|_2^2}{\|A\|_F^2}$ and converges to the least-norm solution $\hat{x}$ of $Ax=b$ in expectation with a linear rate
\begin{equation}\label{equation:3}
\mathbb{E}\|x_{k+1}-\hat{x}\|_2^2\leq (1-\frac{1}{\kappa^2})\mathbb{E}\|x_k-\hat{x}\|_2^2
\end{equation}
where $\kappa=\frac{\|A\|_F}{\sigma_{\min}(A)}$, $\|\cdot\|_F$ is the Frobenius norm,  and $\sigma_{\min}(A)$ is the non-zero smallest singular value of $A$. In order to find sparse solutions of $Ax=b$ instead of the least-norm solutions, the randomized sparse Kaczmarz method (RaSK) was introduced in \citep{2014The,lorenz2014sparse,Petra2015Randomized} with the following scheme
\begin{subequations}
	\begin{align}
	x^*_{k+1}&=x^*_k-\frac{\langle a_i,x_k\rangle-b_i}{\|a_i\|_2^2}a_i,\\
	x_{k+1}&=S_{\lambda}(x^*_{k+1}),
	\end{align}	
	\label{eqation:4}
\end{subequations}
where the index $i$ is sampled by the same rule of RK and $S_{\lambda}(x):=\max\{\lvert x\rvert-\lambda,0\}\cdot \textmd{sign}(x)$ is the soft shrinkage operator with parameter $\lambda\geq 0$. 
It is not hard to see if the parameter $\lambda$ is equal to zero, then the operator $S$ reduces to the identity operator and hence RaSK reduces to RK. 
In this sense, RaSK generalizes RK. Theoretically, it was shown in \cite{2019Linear} that the sequence $\{x_{k}\}$ generated by RaSK converges linearly in expectation to the unique solution $\hat{x}$ of the augmented basis pursuit problem \citep{F2012Exact,2001Atomic,2010Sparse}
\begin{equation}\label{equation:5}
\min_{x\in \mathbb{R}^n}\lambda\|x\|_{1}+\frac{1}{2}\|x\|_2^2,~ \st~ Ax=b.
\end{equation}
Let $supp(\hat{x})=\{j\in \{1,...,n\}\lvert\hat{x}_j\neq 0\}$, and $A_J$ represents the columns of A indexed by $J$, define
\begin{equation}
\label{sigma:1}
\widetilde{\sigma}_{\min}(A)=\min\{\sigma_{\min}(A_J)\lvert J\subset \{1,...,n\},A_J\neq0\}.
\end{equation}
Let $\widetilde{\kappa}=\frac{\|A\|_F}{\widetilde{\sigma}_{\min}(A)},$
and when $b\neq 0$ we have $\hat{x}\neq 0$, define
\begin{equation}
\label{hat:x}
\lvert\hat{x}\rvert_{\min}=\min\{\lvert\hat{x}_j\rvert \lvert j\in supp(\hat{x})\}>0.
\end{equation}
With above notions, the linear convergence rate of RaSK was given in \cite{2019Linear}
\begin{equation}\label{equation:6}
\mathbb{E}\|x_{k+1}-\hat{x}\|_2^2\leq (1- \frac{1}{2}\cdot\frac{1}{\widetilde{\kappa}^2}\cdot \frac{\lvert\hat{x}\rvert_{\min}}{\lvert\hat{x}\rvert_{\min}+2\lambda})\mathbb{E}\|x_k-\hat{x}\|_2^2.
\end{equation}
\par
In view of the Kaczmarz methods use
only a single row of matrix in each iterate,  
they are widely applied to slove large-scale linear systems in various cases, for example, tensor recovery  \citep{chen2021regularized,du2021randomized}, compressed sensing \citep{lorenz2014sparse}, phase retrieval \citep{tan2019phase} and so on. Due to the wide range of applications of Kaczmarz methods, more advanced sampling rules have to be considered such as methods in \citep{2013Randomized,2020A,2021sparse,2021adaptively,2021A,lrr2022} in order to improve the convergence rate of RK or RaSK.
For example, in our recent paper \citep{2021sparse}, we introduced a sparse sampling Kaczmarz-Motzkin method which essentially combined the random and greedy ideas \citep{1992New} together. Very recently, a weighted sampling rule, which selects the $i$-th row with likelihood proportional to $\lvert\langle a_i,x_k\rangle-b_i\rvert^p/\|a_i\|_2^p$ with $0<p<\infty$, was proposed by \citep{2020A,2021A} respectively to accelerate the RK method, and more importantly the new rule can be used to explain why the greedy idea (also called maximal correction method) works well. As a natural question, can we adopt the proposed weighted sampling rule to speedup the RaSK? In this study, we answer this question in an affirmative way.
Actually, we will show that the RaSK equipped with the weighted sampling rule (WRaSK) converges linearly in expectation in the sense that
\begin{equation}\label{equation:9}
\mathbb{E}\|x_{k}-\hat{x}\|_2^2\leq
\left(1-\frac{1}{2} \widetilde{\sigma}^2_{min}(A)\frac{\lvert\hat{x}\rvert_{min}}{\lvert\hat{x}\rvert_{min}+2\lambda}
\mathop{\inf}\limits_{z\neq 0}\frac{\|Az\|_{l_{p+2}}^{p+2}}{\|Az\|_{l_{p}}^{p}\|Az\|_2^{2}}\right)^{k}(2\lambda \|\hat{x}\|_1+\|\hat{x}\|_2^2),
\end{equation}
where $A$ is normalized to $\|a_i\|_2=1$ for each row and $\hat{x}$ is the unique solution of \eqref{equation:5}. By showing
\begin{equation}\label{equation:10}
\mathop{\inf}\limits_{z\neq 0}\frac{\|Az\|_{l_{p+2}}^{p+2}}{\|Az\|_{l_{p}}^{p}\|Az\|_2^{2}}\geq \frac{1}{m}=\frac{1}{\|A\|_F^2},
\end{equation}
we can conclude that the new linear convergence rate is at least as efficient as RaSK; more details will be presented in our Theorem \ref{th1}. On the other hand, WRaSK reduces to RaSK when $p=0$ and approximates the maximal correction variant of RaSK as $p\rightarrow\infty$; the latter also explains why our previously proposed sparse sampling Kaczmarz-Motzkin method in \cite{2021sparse} works well.
Numerically, we demonstrate the superiority of WRaSK via a group of experiments.
\par The paper is organized as follows. In Section \ref{sec:2} we recall the basic knowledge about Bregman distance and Bregman projection. In Section \ref{sec:3} we propose the WRaSK method and prove its linear convergence rates in both noiseless and noisy cases. 
Some detailed remarks about WRaSK are discussed in Section \ref{sec:4}.
In Section \ref{sec:5} we apply some experiments to verify the superiority of WRaSK. Section \ref{sec:6} is our conclusion.

\section{Preliminaries}\label{sec:2}
First, we recall some basic knowledge about convex analysis \citep{2004Introductory,2004Variational}.
\subsection{Convex analysis tool}
\begin{definition}
	Let $f:\mathbb{R}^n\rightarrow \mathbb{R}$ be a convex function. The subdifferential of $f$ at $x\in \mathbb{R}^n$ is
	$$\partial f(x):=\{x^*\in \mathbb{R}^n\lvert f(y)\geq f(x)+\langle x^*,y-x\rangle, \forall y\in \mathbb{R}^n\}.$$
\end{definition}
\noindent
If the convex function $f$ is assumed to be differentiable, then $\partial f(x)=\{\nabla f(x)\}.$

\begin{definition}
	We say that a convex function $f:\mathbb{R}^n\rightarrow \mathbb{R}$ is $\alpha$-strongly convex if there exits $\alpha>0$ such that $\forall x,y\in \mathbb{R}^n$ and $x^*\in \partial f(x)$, we have
	$$f(y)\geq f(x)+\langle x^*,y-x\rangle+\frac{\alpha}{2}\|y-x\|^2.$$
	The conjugate function is $f^*(x^*)=\sup_{x\in\mathbb{R}^{n}}\{\langle x^*,x\rangle-f(x)\}$.
\end{definition}
\noindent
Given a convex function $f$. Then, $f(\cdot)+\frac{\lambda }{2}\|\cdot\|_2^2$ must be $\lambda$-strongly convex if $\lambda>0$.

\begin{example}[\cite{2019Linear}]
	\label{example:1}
	For the augmented $\ell_1$-norm $f(x)=\lambda\|x\|_1+\frac{1}{2}\|x\|_2^2$, its subdifferential is given by
	$\partial f(x)=\{x+\lambda\cdot s\lvert s_i=sign(x_i),x_i\neq 0\ and\ s_i\in[-1,1],x_i=0\}.$
	Its conjugate function is
	$f^*(x)=\frac{1}{2}\|S_{\lambda}(x)\|^2$, and moreover $\nabla f^*(x)=S_{\lambda}(x).$
\end{example}

\subsection{The Bregman distance}
\begin{definition}[\cite{2014The}]
	Let $f:\mathbb{R}^n\rightarrow \mathbb{R}$ be a convex function. The Bregman distance between $x,y \in \mathbb{R}^n$ with respect to f and $x^*\in \partial f(x)$ is defined as
	$$D_f^{x^*}(x,y):=f(y)-f(x)-\langle x^*,y-x\rangle.$$
\end{definition}

\begin{example}
	(a) When $f(x)=\frac{1}{2}\|x\|_2^2$, we have $\partial f(x)=\{\nabla f(x)\}\ and\ D_f^{x^*}(x,y)=\frac{1}{2}\|y-x\|_2^2.$\\
	(b) When $f(x)=\lambda \|x\|_1+\frac{1}{2}\|x\|_2^2$, we have $x^*=x+\lambda\cdot s\in\partial f(x)$ and  $D_f^{x^*}(x,y)=\frac{1}{2}\|y-x\|_2^2+\lambda(\|y\|_1-\langle s,y\rangle),$
	where $s$ is defined in Example \ref{example:1}.
\end{example}
The following result will be used in the forthcoming convergence analysis of WRaSK.
\begin{lemma}[\cite{2019Linear}]\label{lemma:2}
	Let $f:\mathbb{R}^n\rightarrow \mathbb{R}$ be $\alpha$-strongly convex. Then for all $x,y\in \mathbb{R}^n$ and $x^*\in \partial f(x), y^*\in \partial f(y)$, we have
	$$
	\frac{\alpha}{2}\|x-y\|_2^2\leq D_f^{x^*}(x,y)\leq \langle x^*-y^*,x-y\rangle\leq \|x^*-y^*\|_2\cdot\|x-y\|_2
	$$
	and hence
	$D_f^{x^*}(x,y)=0 \Leftrightarrow x=y.$
\end{lemma}

\subsection{The Bregman projection}
\begin{definition}[\cite{2014The}]
	Let $f:\mathbb{R}^n\rightarrow \mathbb{R}$ be $\alpha$-strongly convex and $C\subset \mathbb{R}^n$ be a nonempty closed convex set. The Bregman projection of x onto C with respect to f and $ x^*\in \partial f(x)$ is the unique point $\Pi_{C}^{x^*}(x)\in C$ such that
	\begin{equation}
	D^{x^*}_f(x,\Pi_{C}^{x^*}(x))=\min_{y\in C}D^{x^*}_f(x,y).\notag
	\end{equation}
\end{definition}
The Bregman projection generalizes the traditional orthogonal projection.
Note that when $f(x)=\frac{1}{2}\|x\|_2^2$, the Bregman projection reduces to the  orthogonal projection, denoted by $P_{C}(x)$ as usual.
The next lemma tells us how to compute the Bregman projection onto affine subspaces.
\begin{lemma}[\cite{2014The}]\label{lemma:2.3}
	Let $f:\mathbb{R}^n\rightarrow \mathbb{R}$ be $\alpha$-strongly convex, $a\in \mathbb{R}^n$, $\beta\in \mathbb{R}.$
	The Bregman projection of $x\in \mathbb{R}^n$ onto the hyperplane $H(a,\beta)=\{x\in \mathbb{R}^n\lvert\langle a,x\rangle=\beta\}$ with $a\neq 0$ is
	$$z:=\Pi_{H(a,\beta)}^{x^*}(x)=\nabla f^*(x^*-a\hat{t}),$$
	where $\hat{t}\in \mathbb{R}$ is a solution of
	$\min_{t\in \mathbb{R}}f^*(x^*-at)+t\beta.$
	Moreover, $z^*=x^*-\hat{t}a$ is an admissible subgradient for $z$ and for all $y\in H(a,\beta)$ we have
	\begin{equation}\label{equation:11}
	D^{z^*}_f(z,y)\leq D^{x^*}_f(x,y)-\frac{\alpha}{2}\frac{(\langle a,x\rangle-\beta)^2}{\|a\|_2^2}.
	\end{equation}
\end{lemma}

\section{Weighted Randomized Sparse Kaczmarz method}\label{sec:3}
In this section, we will introduce the WRaSK method to solve the augmented basis pursuit problem (\ref{equation:5}). We assume that $b\neq0$ belongs to $\mathcal{R}(A)$ so that (\ref{equation:5}) has a unique solution $\hat{x}\neq 0$. Without loss of generality, we normalize each row of the matrix $A$ to $\|a_i\|_2=1$.

%Moreover, we will develop the convergence analysis of WRaSK, comparing WRaSK with RaSK in view of convergence rate and time complexity and analyzing the effect of parameter $p$ on the convergence rate of WRaSK.

\subsection{The sampling rule}
Recall that we aim to adopt the weighted sampling rule to RaSK; so let us first introduce what this rule is.
Different from the uniform sampling in RK and RaSK (note that $\|a_i\|_2=1$), the weighted sampling rule gives the rows with large residuals a greater probability. To achieve it mathematically, we start with a given point $x_k$ and then compute the residuals $\lvert\langle a_i,x_k\rangle -b_i\rvert,i=1,...,m$. Since a small residual $\lvert\langle a_i,x_k\rangle -b_i\rvert$ means that $x_k$ approximately solves the equation $\langle a_i,x\rangle =b_i$, so we should try to correct the equations with large residuals. The weighted sampling rule selects the $i$-th row with likelihood proportional to $\lvert\langle a_i,x_k\rangle -b_i\rvert^p$ with $0<p<\infty$. Hence, the rows with greater residual are more possible to be selected. It should be noted that the weighted sampling rule shares the similar greedy idea with the maximal correction method, which selects row with the largest residual in a determined rather than random way.

With the discussion above, WRaSK can be simply obtained by using the scheme \eqref{eqation:4} with the index $i$ being chosen with probability
$p_i=\frac{\lvert\langle a_i,x_k\rangle -b_i\rvert^p}{\|Ax_k-b\|_{l_p}^p}.$

\subsection{The WRaSK method}
In this part, we formally introduce the weighted randomized sparse Kaczmarz method (WRaSK), which combines the randomized sparse Kaczmarz method and the weighted sampling method. It follows that WRaSK inherits their advantages. Note that the iterate $x_k$ generated by WRaSK is projected to the selected hyperplane by using Bregman projection instead of orthogonal projection, which applying the superiority of Bregman distance to obtain sparse solutions of linear systems by using the augmented $l_1$ norm function.
The WRaSK method has two cases, that is, inexact step and exact step, we abbreviate as WRaSK and EWRaSK respectively. In the case of inexact step, $t_k=\langle a_{i_k},x_k\rangle-b_{i_k}$, which  can be seen as a relaxation of Bregman projection. On the other hand, in the case of exact step, we need to solve the piecewise quadratic function minimization problem
\begin{equation}
\label{eq3:1}
t_k=\arg\min_{t\in \mathbb{R}}\frac{1}{2}\|S_{\lambda}(x^*-a_{i_k}t)\|_2^2+tb_{i_k}.
\end{equation}
\par We remark that the computational complexity of WRaSK is less than that of EWRaSK due to the different step-size selections; this point will be further explained in Section \ref{sec:4.1}.
It can be served as the stopping criterion that the maximum iterations step or allowable tolerance error are reached.
\begin{algorithm}
	\caption{Weighted Randomized Sparse Kaczmarz method (WRaSK)}\label{al1}
	\begin{algorithmic}[1]
		\State \textbf{Input}: $x_{0}=x_{0}^*=0\in \mathbb{R}^{n}$,
		$A \in \mathbb{R}^{m \times n}, b \in \mathbb{R}^{m},$ and parameter $p$ 
		\State \textbf{Ouput}: solution of $\min_{x\in \mathbb{R}^n}\lambda\|x\|_{1}+\frac{1}{2}\|x\|_2^2\ ~\st ~Ax=b$
		\State normalize $A$ by row
		\State \textbf{for} $k=0,1,2, \ldots$ do
		\State $\quad$$\quad$choose an index $i_k$ at random with probability 
		$p_i=\frac{ \lvert \langle a_i,x_k\rangle -b_i\rvert ^p}{\|Ax_k-b\|_{l_p}^p}$
		\State $\quad$$\quad$\textbf{switch} Type of step:
		\State $\quad$$\quad$\textbf{case1: inexact step}
		\State $\quad$$\quad$$\quad$$\quad$compute $t_k=\langle a_{i_k},x_k\rangle-b_{i_k}$
		\State $\quad$$\quad$\textbf{case2: exact step}
		\State $\quad$$\quad$$\quad$$\quad$compute $t_k=\arg\min_{t\in \mathbb{R}}f^*(x_k^*-ta_{i_k})+tb_{i_k}$
		\State $\quad$$\quad$\textbf{endSwitch}
		\State $\quad$$\quad$update $x_{k+1}^*=x_k^*-t_ka_{i_k}$
		\State $\quad$$\quad$update $x_{k+1}=S_{\lambda}(x_{k+1}^*)$
		\State $\quad$$\quad$increment $k=k+1$
		\State \textbf{until} a stopping criterion is satisfied
	\end{algorithmic}
\end{algorithm}

\begin{remark}
	Consider two extreme cases of parameter $p$ in the weighted sampling rule.\\
	(a) In the case of $p=0$, we have $p_i=\frac{1}{m},$ which is just the uniform sampling, and hence WRaSK reduces to RaSK under the condition that matrix $A$ is normalized by row. \\
	(b) As $p\rightarrow\infty$, we have
	$i_k=\arg\max_{i=1,...,m}\{\lvert\langle a_i,x_k\rangle -b_i\rvert\},$
	which is the maximal correction method in \cite{1992New}. In this case, WRaSK is a sparse variant of partially randomized Kaczmarz (PRK) in \cite{2020A}.  
\end{remark}

\subsection{Convergence analysis of WRaSK methods}
The following lemma establishes a relationship between $D_f^{x^*}(x,\hat{x})$ and $\|Ax-b\|_2^2$.
\begin{lemma}[\cite{2019Linear}]\label{lemma:3}
	Let $\widetilde{\sigma}_{min}(A)$ and $\lvert\hat{x}\rvert_{min}$ be given in (\ref{sigma:1}), (\ref{hat:x}) respectively.
	Then for all $x\in \mathbb{R}^n$ with $\partial f(x)\cap \mathcal{R}(A^T)\neq \emptyset$
	and all $x^*=A^Ty\in \partial f(x)\cap \mathcal{R}(A^T)$ we have
	$$D_f^{x^*}(x,\hat{x})\leq \frac{1}{\widetilde{\sigma}^2_{min}(A)}\cdot\frac{\lvert\hat{x}\rvert_{min}+2\lambda}{\lvert\hat{x}\rvert_{min}}\cdot\|Ax-b\|_2^2.$$
\end{lemma}
Now, we proceed to show the convergence of WRaSK in noiseless and noisy cases respectively.
\begin{theorem}[Noiseless case]\label{th1}
	Let $0< p< \infty$ and suppose that $A$ is normalized to $\|a_i\|_2=1$ for each row. The sequences $\{x_k\}$ generated by Algorithm 1 converge linearly in expectation to the unique solution $\hat{x}$ of the regularized basis pursuit problem (\ref{equation:5}) in the sense that
	$$\mathbb{E}(D_f^{x_{k+1}^*}(x_{k+1},\hat{x}))\leq
	q\cdot\mathbb{E}(D_f^{x_{k}^*}(x_{k},\hat{x})),
	$$
	and
	$$\mathbb{E}\|x_{k}-\hat{x}\|_2\leq q^{\frac{k}{2}}\sqrt{2\lambda \|\hat{x}\|_1+\|\hat{x}\|_2^2},$$
	where the convergence factor is $q=1-\frac{1}{2}\cdot\widetilde{\sigma}_{\min}^2(A)\cdot\frac{\lvert\hat{x}\rvert_{\min}}{\lvert\hat{x}\rvert_{\min}+2\lambda}\cdot \mathop{\inf}\limits_{z\neq 0}\frac{\|Az\|_{l_{p+2}}^{p+2}}{\|Az\|_{l_{p}}^{p}\|Az\|_2^{2}}$, and $z=x-\hat{x}$.\\
	Moreover, WRaSK is at least as efficient as RaSK due to
	$$\mathop{\inf}\limits_{z\neq0}\frac{\|Az\|_{l_{p+2}}^{p+2}}{\|Az\|_{l_{p}}^{p}\|Az\|_2^{2}}\geq \frac{1}{m},$$
	where the minimum $\frac{1}{m}$ is attainable iff $Az$ is the constant multiple of unit vector.
\end{theorem}
Due to the proof of Theorem \ref{th1} is  involved, it will be given in Appendix \ref{appendix:1}.

\begin{remark}
	We make the following remarks:\\
	(a) WRaSK is at least as efficient as RaSK, independent on the value of $p$.\\
	(b) The convergence factor $q=1-\frac{1}{2}\cdot\widetilde{\sigma}_{\min}^2(A)\cdot\frac{\lvert\hat{x}\rvert_{\min}}{\lvert\hat{x}\rvert_{\min}+2\lambda}\cdot \mathop{\inf}\limits_{z\neq 0}\frac{\|Az\|_{l_{p+2}}^{p+2}}{\|Az\|_{l_{p}}^{p}\|Az\|_2^{2}}$ depends on $\hat{x}$,$\lambda$, $\widetilde{\sigma}_{\min}(A)$ and $p$. If A has full column rank, then $\widetilde{\sigma}_{\min}(A)=\sigma_{\min}(A)$. If $\lambda=0$, then the reliance on $\hat{x}$ disappears.\\
	(c) When $\lvert a_i^Tx_k-b_i\rvert,i=1,...,m$ are equal, the sampling rules of WRaSK and RaSK are the same, in which case the convergence rate of WRaSK is same to that of RaSK.
\end{remark}

In the following, we take the noisy case into account.
\begin{theorem}[Noisy case]\label{th2}
	Assume that a noisy observed data $b^{\delta}\in \mathbb{R}^m$ with $\|b^{\delta}-b\|_{l_{p+2}}\leq\delta$ is given, where $b=A\hat{x}$. Let the sequence $\{x_k\}$ be generated by WRaSK or EWRaSK with $b$ replaced by $b^\delta$.
	Then with the same convergence factor $q$ in Theorem 3.1 as in the noiseless case, we have
	(a) for the WRaSK method:
	$$
	\begin{aligned}
	\mathbb{E}[\|x_k-\hat{x}\|_2]
	&\leq
	q^{\frac{k}{2}}\sqrt{2\lambda \|\hat{x}\|_1+\|\hat{x}\|_2^2}+\delta\sqrt{\frac{c^pq}{1-q}}.
	\end{aligned}
	$$
	(b) for the EWRaSK method:
	$$
	\begin{aligned}
	\mathbb{E}[\|x_k-\hat{x}\|_2]
	&\leq
	q^{\frac{k}{2}}\sqrt{2\lambda \|\hat{x}\|_1+\|\hat{x}\|_2^2}+\delta\sqrt{(1+\frac{4\lambda\|A\|_{1,p+2}}{\delta})\cdot \frac{c^pq}{1-q}}.
	\end{aligned}
	$$
	where constant $c\in\mathbb{R}$ characterizes the equivalence of vector norms $\|\cdot\|_p$ and $\|\cdot\|_{p+2}$.
\end{theorem}
The proof of Theorem \ref{th2} will be given in detail in the Appendix \ref{appendix:2}.

\section{Some remarks about WRaSK}\label{sec:4}
In this section, we would make some noteworthy comments about WRaSK. First, WRaSK is at least efficient as RaSK in terms of convergence rate, and hence we should compare the time complexity of them. 
Second, exploring the effect of parameter $p$ on WRaSK is an interesting problem.
Finally, given that WRaSK needs to compute all residuals in each iterate, 
a partially weighted randomized sparse Kaczmarz is proposed to face with the disadvantage of WRaSK.

\subsection{Time complexity of RaSK and WRaSK}\label{sec:4.1}
The main difference between WRaSK and RaSK is the selection of sampling rules. In this part, we compare them in terms of time complexity. The main time complexity of them in each iteration come from sampling indices and updating the iterative sequences $x_{k}$.
Note that the exact step, obtained by computing the optimal solution of (\ref{eq3:1}), needs $O(n\log(n))$-sorting procedure.
%The sampling complexity of RK or RaSK is not $O(1)$. It is a non-adaptive sampling rule, the probability $p_i$ can be calculated in advance. 
According to Table \ref{ta1}, we conclude that the time complexity of WRaSK and RaSK are in the same level.
Moreover, we have that the time complexity of WRaSK is at the same level as EWRaSK when $m>log(n)$, otherwise, the EWRaSK will cost more time than WRaSK in each iterate.

\begin{table}[h]
	\begin{center}
		\begin{minipage}{\textwidth}
			\caption{Comparisons of Kaczmarz-type methods about time complexity}\label{ta1}
			\begin{tabular*}{\textwidth}{@{\extracolsep{\fill}}lccc@{\extracolsep{\fill}}}
				\toprule%
				Method & Sampling Rule $p_i$& Convergence Rate Bound & Time Complexity\\
				\midrule
				%RK & $p_i=\frac{\|a_i\|_2^2}{\|A\|_F^2}$ & $ 1-\frac{1}{\kappa^2}$ & $O(mn)$\\
				RaSK & $\frac{\|a_i\|_2^2}{\|A\|_F^2}$ &   $1-\frac{1}{2}\cdot  \frac{\lvert\hat{x}\rvert_{min}}{\lvert\hat{x}\rvert_{min}+2\lambda}\cdot\frac{1}{\widetilde{\kappa}^2}$ & $O(mn)$\\
				ERaSK & $\frac{\|a_i\|_2^2}{\|A\|_F^2}$ &   $1-\frac{1}{2}\cdot  \frac{\lvert\hat{x}\rvert_{min}}{\lvert\hat{x}\rvert_{min}+2\lambda}\cdot\frac{1}{\widetilde{\kappa}^2}$ & $O(mn+n\log(n))$\\
				%WRK & $p_i=\frac{|\langle a_i,x_k\rangle -b_i|^p}{\|Ax_k-b\|_{l_p}^p}$ & $1-\mathop{inf}\limits _{z\neq0}\frac{\|Az\|_{l_{p+2}}^{p+2}}{\|Az\|_{l_p}^p\|z\|_2^2}$ & $O(mn+mp)$\\
				WRaSK & $\frac{\lvert\langle a_i,x_k\rangle -b_i\rvert^p}{\|Ax_k-b\|_{l_p}^p}$ & $1-\frac{1}{2}\frac{\lvert\hat{x}\rvert_{min}\cdot\widetilde{\sigma}_{min}^2(A)}{\lvert\hat{x}\rvert_{min}+2\lambda} \mathop{inf}\limits_{z\neq0}\frac{\|Az\|_{l_{p+2}}^{p+2}}{\|Az\|_{l_{p}}^{p}\|Az\|_2^{2}}$ & $O(mn)$ \\
				EWRaSK & $\frac{\lvert\langle a_i,x_k\rangle -b_i\rvert^p}{\|Ax_k-b\|_{l_p}^p}$ &
				$1-\frac{1}{2}\frac{\lvert\hat{x}\rvert_{min}\cdot\widetilde{\sigma}_{min}^2(A)}{\lvert\hat{x}\rvert_{min}+2\lambda} \mathop{inf}\limits_{z\neq0}\frac{\|Az\|_{l_{p+2}}^{p+2}}{\|Az\|_{l_{p}}^{p}\|Az\|_2^{2}}$ & $O(mn+n\log(n))$ \\
				\botrule	
			\end{tabular*}
		\end{minipage}
	\end{center}
\end{table}

\subsection{The effect of parameter $p$ on convergence rate}
Note that the convergence rate of WRaSK is decided by
$$\frac{\|Az\|_{l_{p+2}}^{p+2}}{\|Az\|_{l_{p}}^{p}\|Az\|_2^{2}},$$
whose monotonicity on $p$ was observed by the author in \cite{2021A}. Here, we provide a detailed proof.
We view it as the function about parameter $p$. Redenote $Az=g$, since the quantity will not change when we multiple $g$ by a positive number, we can assume without loss of generality that each non-zero coordinate satisfying $\lvert g_i\rvert\geq1$. Ignoring the constant value $\|g\|_2^2$ and denoting $\log\lvert g_i\rvert=d_i\geq0$, we can reformulate
$$\frac{\|Az\|_{l_{p+2}}^{p+2}}{\|Az\|_{l_{p}}^{p}}=
\frac{\sum_{i=1}^{n}e^{d_i(p+2)}}{\sum_{i=1}^{n}e^{d_ip}}.$$
Now we prove that the function increases as the $p$ increases.
\begin{lemma}
	\label{lemma:4.1}
	If
	$f(x)=\frac{\sum_{i=1}^{n}e^{d_i(x+2)}}{\sum_{i=1}^{n}e^{d_ix}}$, where $d_i\geq0$, $0<x<\infty,$ then $f(x)$ is a monotonic increasing function.
\end{lemma}
The proof will be given in Appendix \ref{appendix:3}.
Lemma \ref{lemma:4.1} verifies that the convergence rate of WRaSK increases as $p$ increases, and hence the increment of parameter $p$ represents the tend of constant improvement of WRaSK.
Moreover, the existence of parameter $p$ of weighted sampling rule connects randomized sampling rule ($p=0$) and greedy sampling rule ($p\rightarrow$ $\infty$). In views of this, there is no need to pursue the optimal parameter $p$. And we will give empirical value of $p$ by developing numerical experiments in Section \ref{sec:5}.

%\begin{remark}
%	(a) Regardless of the value of $p$, WRaSK will perform at least as well as RaSK. In more detail, the greater value of parameter $p$, the faster convergence rate of WRaSK. And this will be verfied by the following experiments.\\
%	(b) The Lemma 3.2 strictly proves why the maximal correction method works well.
%\end{remark}
\subsection{The partially weighted randomized sparse Kaczmarz method}
In our proposed WRaSK, we have to use all information about the linear system to construct the sampling rule. Apparently, the disadvantage of WRaSK is unfavorable for large-scale linear systems, while Kaczmarz has advantage on small storage space and computation in each iterate. Hence, we should consider how to deal with this problem by reducing the number of required residuals in each iterate. 

Guided by the idea of selecting larger residuals, a randomized Kaczmarz with partially weighted selection step was proposed in
\cite{gross2021note}. 
It constantly compares the residuals of the sampling indices and leaves a larger residual.
In order to overcome the weakness of WRaSK, we provide a possible solution, which is
a combination of randomized sparse Kaczmarz and partially weighted selection. We propose the partially weighted randomized sparse Kaczmarz method (PWRaSK).
The pseudocode is as follows.

Denote the absolute residual of $i$-th index in $k$-th iterate as $r_k(i),1\leq i\leq m$. Note that the step 5-8 in PWRaSK, which is to find a relatively large residual. The sampling method samples indices consistently, and then compares the new residual with the candidate residual. If the new is smaller, the candidate is used as iterative index, otherwise, the new is used as the candidate to continue.
Similar to the convergence rate of the algorithm proposed by \cite{gross2021note}, the convergence rate of PWRaSK can be generalized from randomized sparse Kaczmarz (RaSK) \citep{2014The,patel2021convergence}. Moreover, it is expected that PWRaSK converges faster than RaSK due to its sampling rule approaches the maximum residual method \citep{1992New} or partially randomized Kaczmarz \citep{2020A}.
\begin{algorithm}
	\caption{Partially Weighted Randomized Sparse Kaczmarz (PWRaSK)}\label{al2}
	\begin{algorithmic}[1]
		\State \textbf{Input}: $x_{0}=x_{0}^*=0\in \mathbb{R}^{n}$,
		$A \in \mathbb{R}^{m \times n}, b \in \mathbb{R}^{m},$ and parameter $p$ 
		\State \textbf{Ouput}: solution of $\min_{x\in \mathbb{R}^n}\lambda\|x\|_{1}+\frac{1}{2}\|x\|_2^2\ ~\st ~Ax=b$
		\State normalize $A$ by row
		\State \textbf{for} $k=0,1,2, \ldots$ do
		\State $\quad$$\quad$set $U=\{1,\cdots,m\}$, select an index $i$ from $U$ as following method:
		\State $\quad$$\quad$select an index $i_1$ from $U$ uniformly
		\State $\quad$$\quad$set $U=U\backslash\{i_1\}$, if $U=\emptyset$ set $i=i_1$ and go to Step 9,
		otherwise, select an index $i_2$ from $U$ uniformly
		\State $\quad$$\quad$if $r_k(i_1)>r_k(i_2)$, set $i=i_1$ and go to Step 9. Otherwise, set $i_1=i_2$ and go to Step 7
		\State $\quad$$\quad$\textbf{switch} Type of step:
		\State $\quad$$\quad$\textbf{case1: inexact step}
		\State $\quad$$\quad$$\quad$$\quad$compute $t_k=\langle a_{i_k},x_k\rangle-b_{i_k}$
		\State $\quad$$\quad$\textbf{case2: exact step}
		\State $\quad$$\quad$$\quad$$\quad$compute $t_k=\arg\min_{t\in \mathbb{R}}f^*(x_k^*-ta_{i_k})+tb_{i_k}$
		\State $\quad$$\quad$\textbf{endswitch}
		\State $\quad$$\quad$update $x_{k+1}^*=x_k^*-t_ka_{i_k}$
		\State $\quad$$\quad$update $x_{k+1}=S_{\lambda}(x_{k+1}^*)$
		\State $\quad$$\quad$increment $k=k+1$
		\State \textbf{until} a stopping criterion is satisfied
	\end{algorithmic}
\end{algorithm}

\section{Numerical experiments}\label{sec:5}
In this section, we will discuss the effectiveness of WRaSK by experiments. In section \ref{sec:5.1}, we analyze the effect of parameter $p$ on WRaSK. In section \ref{sec:5.2}, we are going to solve linear systems with the coefficient matrix $A\in \mathbb{R}^{m\times n}$ being generated by MATLAB function 'randn' or chosen from SuiteSparse Matrix Collection in \cite{davis2011university}.

Without loss of generality, we normalize A by row. Moreover, $\hat{x}\in \mathbb{R}^n$ is created by using the MATLAB function 'sparserandn' and randomly choosing the nonzero location by sparsity.
The exact data $b\in \mathbb{R}^m$ is calculated by $A\hat{x}=b$ in noiseless case, while we use 'randn' to generate the noise $r$ and then we obtain noisy data $b^{\delta} =b+r$.
Our experiments are initialized with $x_0=x_0^*=0$, and carry out 60 trials to ensure accuracy. The main outputs are relative residual and relative error,  defined respectively by
$$Residaul =\frac{\|Ax-b\|_2}{\|b\|_2}, Error=\frac{\|x-\hat{x}\|_2}{\|\hat{x}\|_2}.$$
All experiments are performed with MATLAB (version R2021b) on a personal computer with 2.80-GHZ
CPU(Intel(R) Core(TM) i7-1165G7), 16-GB memory, and Windows operating system(Windows 10). 
\subsection{Parameter tuning}\label{sec:5.1}
In this test, numerical experiments are developed to verify Lemma \ref{lemma:4.1}. Set $n=200$, $m=400$, $sparsity =25$, and then select $p$ from $1,\frac{m}{80},\frac{m}{40},\infty$ in WRaSK. Moreover,
matrix $A\in \mathbb{R}^{m\times n}$ is generated by MATLAB function 'randn'.
%Note that WRaSK simply becomes the maximal correction method when $p\rightarrow\infty$, while WRaSK with $p=0$ is equivalent to RaSK under the condition that all rows are normalized.
To explore the performance of WRaSK with different values of $p$, we run 60 times to compute the median of relative residuals and errors from RaSK ($p=0$), WRaSK with different parameters $p$ and maximal correction method ($p\rightarrow\infty$).
\par According to the Fig. 1, we observe that WRaSK performs better as $p$ increases and finally approximates the maximal correction method, which confirms Lemma \ref{lemma:4.1}. 
%Actually, WRaSK provides a connection between randomized Kaczmarz method ($p=0$) and greedy Kaczmarz method ($p\rightarrow \infty$) by changing parameter $p$. 
Since we do not pay attention to the optimal value of parameter $p$, and $p=\frac{m}{40}$ is good enough according to the experiment, we take $p=\frac{m}{40}$ for WRaSK as empirical value in the following experiments.

\begin{figure}[h]\label{fig:1}
	\centering
	\subfigure[The relative residuals]{
		\includegraphics[width=0.43\linewidth,height=0.35\textwidth]{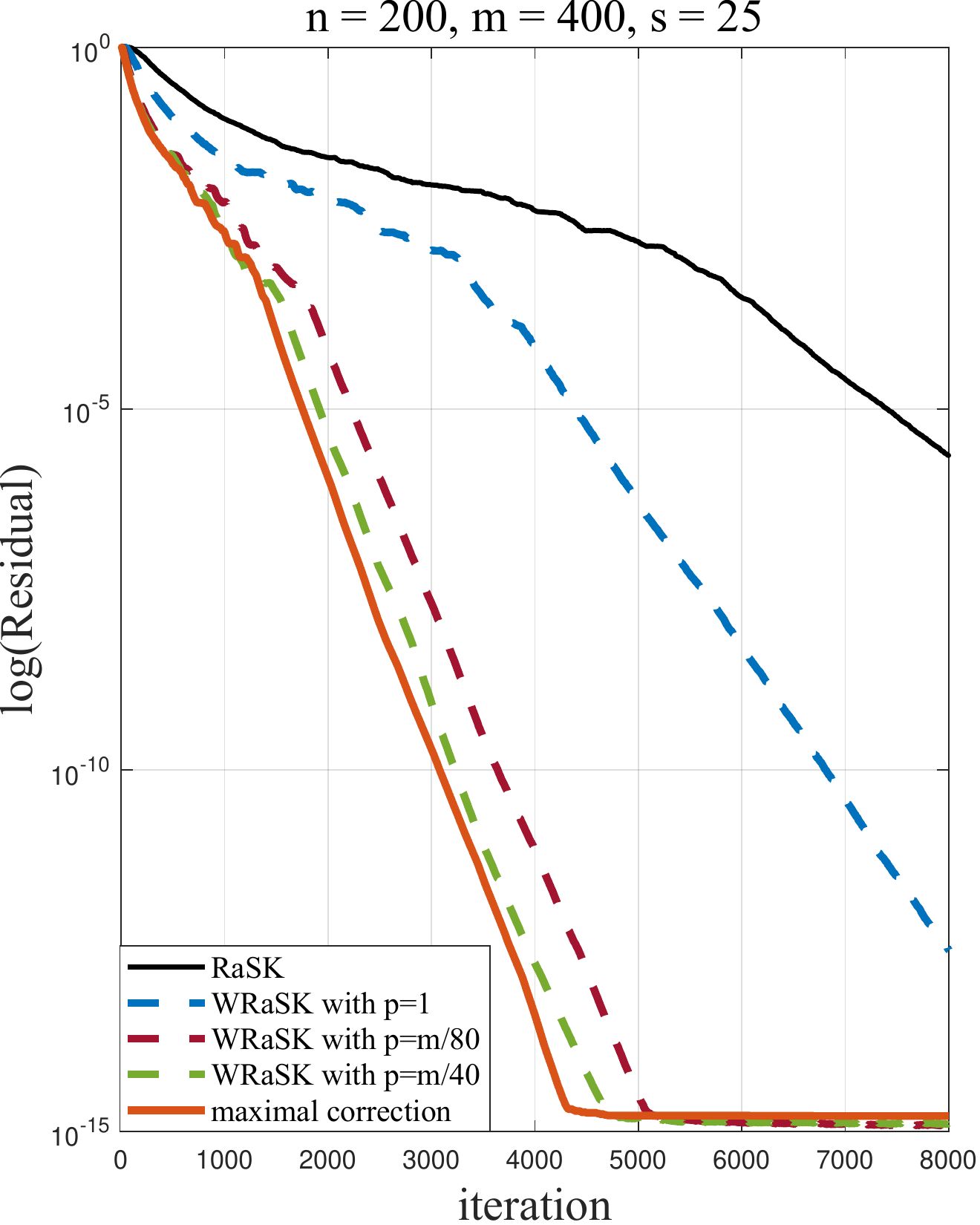}}
	\subfigure[The relative errors]{
		\includegraphics[width=0.43\linewidth,height=0.35\textwidth]{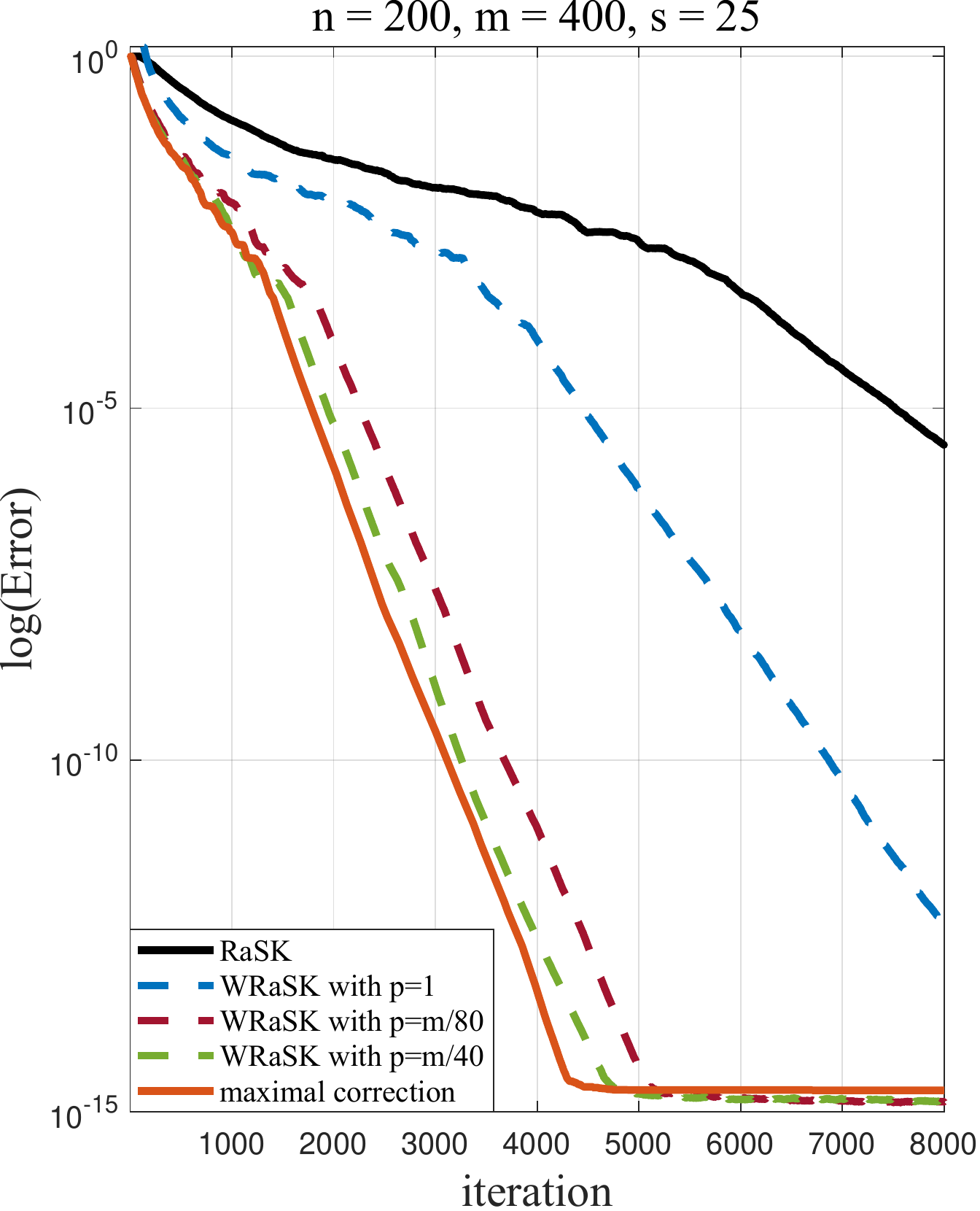}}
	\caption{
		Compare the residuals and errors of WRaSK with different parameters $p$, i.e. RaSK (black), $p$=1 (blue), $p$=m/80 (red), $p$=m/40 (green), and $p\rightarrow\infty $ (orange),
		$n=200,m=400,sparsity=25,\lambda=1$. Left: plots of realtive residual $\|Ax-b\|_2/\|b\|_2$, right: plots of realtive error $\|x-\hat{x}\|_2/\|x\|_2$.	
		The line shows median over 60 trials.
	}
\end{figure}

\subsection{Comparisons of different Kaczmarz variants}\label{sec:5.2}
Here, we explore the performance of WRaSK in terms of  
the number of iteration steps (IT) and the computing time in seconds (CPU). We perform 60 numerical experiments and record IT and CPU when relative error reach given accuracy or the maximum number of iteration steps is achieved. Then take the median of CPU and IT separately in 60 records.

\subsubsection{Simulated data} 
To compare RaSK, WRK and WRaSK, we construct coefficient matrix $A$ by MATLAB 'randn'. 
%For generality, we compare the IT and CPU of these methods in four settings of linear systems respectively.
Set $p=\frac{m}{40}, sparsity=20,\lambda=1$.
Moreover, the iterative processing determinates once the relative error satisfies $Error<10^{-3}$ in noiseless case and $Error<10^{-1}$ in noisy case, or the maximum number of iteration steps reaches 200000.  
If the number of iteration steps achieves 200000, we denote it as "-".
The experimental results are recorded in Table \ref{table:2}.

\begin{table}[h]
	\begin{center}
		\begin{minipage}{\textwidth}
	\caption{IT and CPU of Kaczmarz-type methods for matrices $A$ generated from 'randn'}
	\label{table:2}
	\begin{tabular*}{\textwidth}{@{\extracolsep{\fill}}lccccc@{\extracolsep{\fill}}}
		\toprule
		\multicolumn{2}{c}{} &  \multicolumn{2}{c}{noiseless} & \multicolumn{2}{c}{noisy }\\
		\hline \multicolumn{2}{c}{$m \times n$} & $2000 \times 500$ & $500 \times 2000$ & $1000 \times 500$ & $500 \times 1000$  \\
		\midrule
        \multirow{2}{*}{ ERaSK } & IT & 5290.1 & 160250 & 3810 & 108980 \\
		\cline { 2 - 6} & CPU & 51.2570 & 5797 & 35.2 & 1110.6 \\
		\midrule
		\multirow{2}{*}{ WRK } & IT & 1689.1 & - & - &  -\\
		\cline { 2 - 6 } & CPU & 26.1041 & - & - &  -\\
		\midrule
		\multirow{2}{*}{ EWRaSK } & IT & 588.4 & 117330 & 390 &  113790\\
		\cline { 2 - 6 } & CPU & 9.2601 & 11217 & 5.6 &  1778\\
		\botrule
	\end{tabular*}
\end{minipage}
\end{center}
\end{table}
According to Table \ref{table:2}, we find that whether there is noise or not, EWRaSK 
has the best performance in overdetermined cases, while RaSK outperforms other methods in underdetermined case. It requires the least CPU time and iteration steps to satisfy the stopping rules.
Note that WRK has the worst performance in noisy or underdetermined case, where it fails to converge.

\subsubsection{Real data}

In this subsection, the performance of WRaSK, PWRaSK and WSSKM will be verified under some real matrices.
The coefficient matrices $A$ come from the SuiteSparse Matrix Collection \citep{davis2011university}, which are originated in different applications such as least squares problem, 2D/3D problem, optimal control problem, and so on. 

In this test, we generate ground truth $\hat{x}$ with $sparsity=20$ and compute $A\hat{x}$ as observed data $b$, and the datails of chosen matrices are presented in Table \ref{table:3}.
We perform ERaSK, WRK and EWRaSK on these problems, and determinate algorithms once relative error reaches $10^{-3}$ or the maximum iterations attains 200000.
The IT and CPU of the median of 60 trails in different settings are listed in Table \ref{table:3}.
\begin{table}[h]
	\begin{center}
		\begin{minipage}{\textwidth}
	\caption{IT and CPU of Kaczmarz-type methods for matrices $A$ from real data}
	\label{table:3}
	\begin{tabular*}{\textwidth}{@{\extracolsep{\fill}}lcccc@{\extracolsep{\fill}}}
		\toprule
		\multicolumn{2}{c}{ Name } & WorldCities  &  Trefethen\_300 & Trefethen\_700 \\
		\midrule
		\multicolumn{2}{c}{ $m\times n$ } & $315\times 100$  & $300\times 300$ & $700\times 700$	
		\\
		\multicolumn{2}{c}{Density} & 23.87\% & 5.20\% & 2.58\%
		\\	
		\multicolumn{2}{c}{Cond(A)} & 66.00  & 1772.69 & 4710.40
		\\
		\midrule
		\multirow{2}{*}{ ERaSK } & IT & 4850.5  & 2256 & 9127.5\\
		\cline { 2 - 5} & CPU & 0.1653  & 0.1165 & 0.7257\\
		\midrule
		\multirow{2}{*}{ WRK } & IT & 5414.5  & 183 & 134.5\\
		\cline { 2 - 5} & CPU & 0.4126  & 0.0122 & 0.0189\\
		\midrule
		\multirow{2}{*}{ EWRaSK } & IT & 708  & 24 & 21\\
		\cline { 2 - 5} & CPU & 0.0541  & 0.0014 & 0.0015\\
		\botrule
	\end{tabular*}
\end{minipage}
\end{center}
\end{table}
According to Table \ref{table:3}, we obtain that EWRaSK performs well in terms of IT and CPU in total. Although the sampling rule of WRaSK costs more time in each iterate, the overall time of WRaSK is less than the others.

\section{Conclusion}\label{sec:6}
In this paper, we propose WRaSK method for finding sparse solutions of consistent linear systems, along with detailed analysis of convergence rate in both noiseless and noisy cases. WRaSK, as a variant of RaSK, combines the advantages of RaSK and the weighted sampling method.
The theoretical results show that WRaSK is at least as efficient as RaSK.  
In order to overcome the cost caculation of WRaSK, we provide a possible solution: PWRaSK, which reduces the calculation of residuals used in each iterate.
Numerical experiments also demonstrate the superiority of WRaSK. 
%We conclude the paper with some possible future research directions.

As future work, we wonder whether it is possible to extend our study to other recently proposed sparse Kaczmarz methods.

\bmhead{Acknowledgements}
We would like to thank Dr. Lionel N. Tondji et.al. for sending us their conference paper \cite{tondjilinear} to kindly remind us that they also independently proposed the weighted randomized sparse Kaczmarz method.
This work was supported by the National Natural Science Foundation of China (No.11971480, No.61977065), the Natural Science Fund of Hunan for Excellent Youth (No.2020JJ3038), and the Fund for NUDT Young Innovator Awards (No.20190105).

\section*{Statements and Declarations}
There is no conflict of interest in the manuscript.
The data used in the manuscript are available in the SuiteSparse Matrix Collection.
All authors contributed to the study conception and design. The first draft of the manuscript was written by Lu Zhang and all authors commented on previous versions of the manuscript. All authors read and approve the final manuscript, and are all aware of the current submission to COAM. 

\begin{appendices}
\section{Proof of Theorem 3.1}\label{appendix:1}
\begin{proof}
	The proof is divided into two parts, we deduce the convergence rate of WRaSK in the first part and compare the convergence rate between WRaSK and RaSK in the second part.
	
	First, we derive the convergence rate of WRaSK. By Theorem 2.8 in \cite{2014The} we know that (\ref{equation:11}) in Lemma \ref{lemma:2.3} holds for both the exact and inexact stepsize. Note that $f$ is 1-strongly convex and $\|a_{i_k}\|_2=1$, it follows that
	\begin{align}
	\label{chapter3:eq1}
	D_f^{x_{k+1}^*}(x_{k+1},\hat{x})
	\leq D_f^{x_{k}^*}(x_{k},\hat{x})-\frac{1}{2}(\langle a_{i_k},x_k\rangle-b_{i_k})^2,
	\end{align}
	we fix the values of the indices $i_0,...,i_{k-1}$ and only consider $i_k$ as a random variable. Taking the conditional expectation on both sides we derive that
	\begin{align}\label{eq:1}
	&~~~\mathbb{E} (D_{f}^{x_{k+1}^*}(x_{k+1},\hat{x})\lvert i_{0},...,i_{k-1})\notag\\
	&\leq
	D_{f}^{x_{k}^*}(x_{k},\hat{x})-\frac{1}{2}\sum_{i=1}^{m}(\langle a_{i},x_{k}\rangle-b_{i})^2\cdot\frac{\lvert\langle a_{i},x_{k}\rangle -b_{i}\rvert^p}{\|Ax_{k}-b\|_{l_p}^p}\notag\\
	&=D_f^{x_{k}^*}(x_{k},\hat{x})-\frac{1}{2}\frac{\|Ax_k-b\|_{l_{p+2}}^{p+2}}{\|Ax_k-b\|_{l_{p}}^{p}\|Ax_k-b\|_2^{2}}\cdot\|Ax_k-b\|_2^{2}\notag\\
	&\leq  D_f^{x_{k}^*}(x_{k},\hat{x})-\frac{1}{2}\mathop{\inf}\limits_{x\neq\hat{x}}\frac{\|Ax-b\|_{l_{p+2}}^{p+2}}{\|Ax-b\|_{l_{p}}^{p}\|Ax-b\|_2^{2}}\cdot\|Ax_k-b\|_{2}^{2}\notag\\
	&\leq
	\left(1-\frac{1}{2} \widetilde{\sigma}^2_{\min}(A)\cdot\frac{\lvert\hat{x}\rvert_{\min}}{\lvert\hat{x}\rvert_{\min}+2\lambda}\cdot
	\mathop{\inf}\limits_{x\neq\hat{x}}\frac{\|Ax-b\|_{l_{p+2}}^{p+2}}{\|Ax-b\|_{l_{p}}^{p}\|Ax-b\|_2^{2}}\right)
	D_{f}^{x_{k}^*}(x_{k},\hat{x}).\notag
	\end{align}
	The last inequality follows by invoking Lemma \ref{lemma:3}.
	Now considering all indices $i_0,...,i_k$ as random variables and taking the full expectation on both sides, we have that
	$$\mathbb{E}(D_f^{x_{k+1}^*}(x_{k+1},\hat{x}))\leq
	\left(1-\frac{1}{2} \widetilde{\sigma}^2_{\min}(A)\cdot\frac{\lvert\hat{x}\rvert_{\min}}{\lvert\hat{x}\rvert_{\min}+2\lambda}\cdot\mathop{\inf}
	\limits_{z\neq 0}\frac{\|Az\|_{l_{p+2}}^{p+2}}{\|Az\|_{l_{p}}^{p}\|Az\|_2^{2}}\right)\mathbb{E}(D_f^{x_{k}^*}(x_{k},\hat{x})),
	$$
	where $z=x-\hat{x}$.
	According to Lemma \ref{lemma:2} and $f$ is 1-strongly convex, we can obtain
	\begin{equation}
	D_f^{x_{k}^*}(x_{k},\hat{x})\geq \frac{1}{2}\|x_{k}-\hat{x}\|_2^2.\notag
	\end{equation}
	Thus, we get
	$$\mathbb{E}\|x_{k}-\hat{x}\|_2\leq
	\left(1-\frac{1}{2} \widetilde{\sigma}^2_{\min}(A)\cdot\frac{\lvert\hat{x}\rvert_{\min}}{\lvert\hat{x}\rvert_{\min}+2\lambda}\cdot\mathop{\inf}
	\limits_{z\neq 0}\frac{\|Az\|_{l_{p+2}}^{p+2}}{\|Az\|_{l_{p}}^{p}\|Az\|_2^{2}}\right)^{\frac{k}{2}}
	\sqrt{2\lambda \|\hat{x}\|_1+\|\hat{x}\|_2^2}.
	$$
	Next we compare the convergence rates between RaSK and WRaSK.  H\"older's inequality
	implies that for any $0\neq x\in \mathbb{R}^m,$
	\begin{equation}\label{holder:1}
	\|x\|_{l_p}^p=\sum_{i=1}^{m}\lvert x_i\rvert^p\leq (\sum_{i=1}^{m}\lvert x_i\rvert^{p+2})^{\frac{p}{p+2}}(\sum_{i=1}^{m}1)^{\frac{2}{p+2}}=\|x\|_{l_{p+2}}^pm^{\frac{2}{p+2}},	
	\end{equation}
	and
	\begin{equation}\label{holder:2}
	\|x\|_{2}^2=\sum_{i=1}^{m}\lvert x_i\rvert^{2}\leq (\sum_{i=1}^{m}\lvert x_i\rvert^{p+2})^{\frac{2}{p+2}}(\sum_{i=1}^{m}1)^{\frac{p}{p+2}}=\|x\|_{l_{p+2}}^2m^{\frac{p}{p+2}}.	
	\end{equation}
	Based on (\ref{holder:1}) and (\ref{holder:2}), for $0\neq Az\in\mathbb{R}^m$ we deduce that
	\begin{equation}
	\label{two:inequality}
	\frac{\|Az\|_{l_{p+2}}^{p+2}}{\|Az\|_{l_{p}}^{p}}
	\geq \frac{\|Az\|_{l_{p+2}}^{2}}{m^{\frac{2}{p+2}}}
	\geq \frac{1}{m}\|Az\|_2^2.
	\end{equation}
	Hence,
	\begin{equation}
	\label{inequality:holder}
	\frac{\|Az\|_{l_{p+2}}^{p+2}}{\|Az\|_{l_{p}}^{p}\|Az\|_2^{2}}\geq \frac{1}{m}.
	\end{equation}
	It follows that
	$$\mathop{\inf}\limits_{z\neq0}\frac{\|Az\|_{l_{p+2}}^{p+2}}{\|Az\|_{l_{p}}^{p}\|Az\|_2^{2}}\geq \frac{1}{m},$$
	with which we further derive that
	$$\frac{1}{2}\cdot \widetilde{\sigma}^2_{\min}(A)\cdot\frac{\lvert\hat{x}\rvert_{\min}}{\lvert\hat{x}\rvert_{\min}+2\lambda}\cdot\mathop{\inf}
	\limits_{z\neq0}\frac{\|Az\|_{l_{p+2}}^{p+2}}{\|Az\|_{l_{p}}^{p}\|Az\|_2^{2}}\geq
	\frac{1}{2}\cdot\frac{1}{m}\cdot \widetilde{\sigma}^2_{\min}(A)\cdot\frac{\lvert\hat{x}\rvert_{\min}}{\lvert\hat{x}\rvert_{\min}+2\lambda}.$$
	Thereby, we conclude that the convergence rate of WRaSK is at least as efficient as RaSK.
	
	As we all known, H\"older's inequality takes the equal sign if and only if one of the two vectors is the constant multiple of the other.
	As for cases of equality, it follows from inequality (\ref{inequality:holder}) is deduced by using H\"older's inequality twice that the equality holds if and only if $Az$ is the constant multiple of unit vector.
	The proof is completed.	
\end{proof}		

\section{Proof of Theorem 3.2}\label{appendix:2}
\begin{proof}
	Making use of the observation in \cite{2010Randomized} that
	\begin{equation}
	\label{obser}
	x_k^{\delta}:=\hat{x}+\frac{b_{i_k}^{\delta}-b_{i_k}}{\|a_{i_k}\|_2^2}a_{i_k} \in H(a_{i_k},b_{i_k}^{\delta}).
	\end{equation}
	Note that $f$ is 1-strongly convex and $\|a_{i_k}\|_2=1$, hence according to Lemma \ref{lemma:2.3} we deduce that
	\begin{equation}
	\label{change:Lemma2:3}
	D_f^{x_{k+1}^*}(x_{k+1},x_k^{\delta})\leq D_f^{x_{k}^*}(x_{k},x_k^{\delta})-\frac{1}{2}(\langle a_{i_k},x_k\rangle-b_{i_k}^{\delta})^2.
	\end{equation}
	Reformulating (\ref{change:Lemma2:3}) by
	(\ref{obser}), we derive that
	\begin{equation}\label{eq:2}
	D_f^{x_{k+1}^*}(x_{k+1},\hat{x})\leq
	D_f^{x_{k}^*}(x_{k},\hat{x})-\frac{1}{2}(\langle a_{i_k},x_k\rangle-b_{i_k}^{\delta})^2+\langle x_{k+1}^*-x_k^*,x_k^{\delta}-\hat{x}\rangle.
	\end{equation}
	(a) In the WRaSK method, we have
	$$x_{k+1}^*-x_k^*=-(\langle a_{i_k},x_k\rangle-b_{i_k}^{\delta})a_{i_k}.$$
	Recall that $x_k^{\delta}-\hat{x}=(b_{i_k}^{\delta}-b_{i_k})a_{i_k}$,
	we get
	\begin{equation}
	\langle x_{k+1}^*-x_k^*,x_k^{\delta}-\hat{x}\rangle
	=(b_{i_k}^{\delta}-b_{i_k})^2
	-(b_{i_k}^{\delta}-b_{i_k})\cdot(\langle a_{i_k},x_k\rangle -b_{i_k}),\label{eq:4}
	\end{equation}
	and
	\begin{align}
	\label{eq:11}
	&~~~~-\frac{1}{2}(\langle a_{i_k},x_k\rangle-b_{i_k}^{\delta})^2\nonumber\\
	&=
	-\frac{1}{2}(\langle a_{i_k},x_k\rangle-b_{i_k}+b_{i_k}-b_{i_k}^{\delta})^2\nonumber\\
	&=-\frac{1}{2}(\langle a_{i_k},x_k\rangle-b_{i_k})^2+
	(b_{i_k}^{\delta}-b_{i_k})\cdot(\langle a_{i_k},x_k\rangle-b_{i_k})
	-\frac{1}{2}(b_{i_k}-b_{i_k}^{\delta})^2.
	\end{align}
	Plugging the reformulations (\ref{eq:4}) and (\ref{eq:11}) into (\ref{eq:2}) we have
	$$D_f^{x_{k+1}^*}(x_{k+1},\hat{x})\leq
	D_f^{x_{k}^*}(x_{k},\hat{x})-\frac{1}{2}(\langle a_{i_k},x_k\rangle-b_{i_k})^2
	+\frac{1}{2}(b_{i_k}-b_{i_k}^{\delta})^2.
	$$
	We fix the values of the indices $i_0,...,i_{k-1}$ and only consider $i_k$ as a random variable. Taking the conditional expectation on both sides we get
	$$
	\begin{aligned}
	&~~~~	\mathbb{E}(D_f^{x_{k+1}^*}(x_{k+1},\hat{x})\lvert i_0,...,i_{k-1})\\
	&\leq D_f^{x_{k}^*}(x_{k},\hat{x})-
	\frac{1}{2}\sum_{i=1}^{m}(\langle a_i,x_k\rangle-b_i)^2\frac{\lvert\langle a_i,x_k\rangle -b_i\rvert^p}{\|Ax_k-b\|_{l_p}^p}+
	\frac{1}{2}\sum_{i=1}^{m}(b_i-b_i^{\delta})^2\frac{\lvert\langle a_i,x_k\rangle -b_i\rvert^p}{\|Ax_k-b\|_{l_p}^p}\notag\\
	&\leq qD_f^{x_{k}^*}(x_{k},\hat{x})+\frac{1}{2}\|b-b^{\delta}\|_{l_{p+2}}^2\cdot\frac{\|Ax_k-b\|_{l_{p+2}}^p}{\|Ax_k-b\|_{l_{p}}^p}.
	\end{aligned}$$
	The last inequality can be deduced by using the conclusion of Theorem \ref{th1} and H\"older's inequality
	\begin{equation}
	\sum_{i=1}^{m}(b_i-b_i^{\delta})^2\cdot\lvert\langle a_i,x_k\rangle -b_i\rvert^p\leq \|b-b^{\delta}\|_{l_{p+2}}^2\cdot\|Ax_k-b\|_{l_{p+2}}^p.\notag
	\end{equation}			
	Now considering all indices $i_0,...,i_k$ as random variables and taking the full expectation on both sides, we can derive that
	$$\begin{aligned}
	\mathbb{E}(D_f^{x_{k+1}^*}(x_{k+1},\hat{x}))
	&\leq q^{k+1}(\lambda \|\hat{x}\|_1+\frac{1}{2}
	\|\hat{x}\|_2^2)+\frac{1}{2}\|b-b^{\delta}\|_{l_{p+2}}^2\sum_{i=0}^{k}q^{k-i}
	\frac{\|Ax_i-b\|_{l_{p+2}}^p}{\|Ax_i-b\|_{l_{p}}^p}.
	\end{aligned}$$
	According to the equivalence of vector norms in $\mathbb{R}^m$, there is a constant $c\in\mathbb{R}$ such that for any vector $z\in\mathbb{R}^m$ we have that
	\begin{equation}
	\|z\|_{l_{p+2}}\leq c\|z\|_{l_{p+2}}.\nonumber
	\end{equation}
	Thus,
	$$\begin{aligned}
	\mathbb{E}(D_f^{x_{k}^*}(x_{k},\hat{x}))
	&\leq q^k(\lambda \|\hat{x}\|_1+\frac{1}{2}
	\|\hat{x}\|_2^2)+\frac{1}{2}\|b-b^{\delta}\|_2^2\cdot \frac{c^pq}{1-q}.
	\end{aligned}$$
	Using $\sqrt{u+v}\leq\sqrt{u}+\sqrt{v}$ and $f$ is 1-strongly convex,  we further deduce that
	$$
	\begin{aligned}
	\mathbb{E}\|x_{k}-\hat{x}\|
	&\leq
	q^{\frac{k}{2}}\sqrt{2\lambda \|\hat{x}\|_1+
		\|\hat{x}\|_2^2}+
	\delta\sqrt{ \frac{c^pq}{1-q}}.\\
	\end{aligned}
	$$
	(b) In the EWRaSK method, according to Example \ref{example:1} we have $x_{k+1}^*=x_k+\lambda\cdot s_k$, where $\|s_k\|_{\infty},\|s_{k+1}\|_{\infty}\leq1$. The exact linesearch guarantees
	$\langle x_{k+1},a_{i_k}\rangle=b_{i_k}^{\delta};$ thus
	\begin{align}
	\label{noise:1}
	\langle x_{k+1}^{*}-x_{k}^{*},x_k^{\delta}-\hat{x}\rangle
	&=\frac{b_{i_{k}}^{\delta}-b_{i_{k}}}{\|a_{i_{k}}\|_2^2}(\langle x_{k+1}-x_{k},a_{i_{k}}\rangle +\lambda\langle s_{k+1}-s_{k},a_{i_k}\rangle),\notag\\
	&\leq \frac{(b_{i_{k}}^{\delta}-b_{i_{k}})^{2}}{\|a_{i_{k}}\|_2^{2}}
	-\frac{(b_{i_{k}}^{\delta}-b_{i_{k}})(\langle a_{i_{k}},x_{k}\rangle-b_{i_{k}})}{\|a_{i_{k}}\|_2^{2}}+
	\frac{2\lambda\lvert b_{i_{k}}^{\delta}-b_{i_{k}}\rvert\cdot\|a_{i_{k}}\|_{1}}{\|a_{i_{k}}\|_{2}^{2}}.
	\end{align}
	Bringing the (\ref{eq:11}) and  (\ref{noise:1}) into (\ref{eq:2}), note that $\|a_{i_k}\|_2=1$ we derive
	\begin{align}
	D_f^{x_{k+1}^*}(x_{k+1},\hat{x})
	&\leq
	D_f^{x_{k}^*}(x_{k},\hat{x})-\frac{1}{2}(\langle a_{i_k},x_k\rangle-b_{i_k})^2
	+\frac{1}{2}(b_{i_k}-b_{i_k}^{\delta})^2+
	2\lambda\lvert b_{i_k}^{\delta}-b_{i_k}\rvert\cdot\|a_{i_k}\|_1.
	\end{align}
	Use H\"older's inequality to reformulate
	\begin{align}
	&2\lambda\sum_{i=1}^{m}\lvert b_{i}^{\delta}-b_{i}\rvert\cdot\|a_{i}\|_1\cdot\frac{\lvert\langle a_i,x_k\rangle -b_i\rvert^p}{\|Ax_{k}-b\|_{l_p}^p}
	&\leq
	2\lambda\|b^{\delta}-b\|_{l_{p+2}}\cdot \|A\|_{1,p+2}\cdot\frac{\|Ax_k-b\|_{l_{p+2}}^p}{\|Ax_k-b\|_{l_p}^p},
	\end{align}
	Similar to (a), we get
	$$
	\begin{aligned}
	\mathbb{E}[\|x_k-\hat{x}\|_2]
	&\leq
	q^{\frac{k}{2}}\sqrt{2\lambda \|\hat{x}\|_1+\|\hat{x}\|_2^2}+\delta\sqrt{(1+\frac{4\lambda\|A\|_{1,p+2}}{\delta})\cdot \frac{c^pq}{1-q}}.
	\end{aligned}
	$$
	The proof is completed.
\end{proof}

\section{Proof of Lemma 4.1}\label{appendix:3}
\begin{proof}
	First, we compute the derivative of the function $f(x)$ as follows
	$$
	\begin{aligned}
	f'(x)
	&=
	\frac{1}
	{(\sum_{i=1}^{n}e^{d_ix})^2}
	[(\sum_{i=1}^{n}d_ie^{2d_i}e^{d_ix})(\sum_{j=1}^{n}e^{d_jx})-
	(\sum_{i=1}^{n}e^{2d_i}e^{d_ix})(\sum_{j=1}^{n}d_je^{d_jx})],\\
	&=
	\frac{1}
	{(\sum_{i=1}^{n}e^{d_ix})^2}
	[\sum_{i=1}^{n}\sum_{j=1}^{n}d_ie^{2d_i}e^{(d_i+d_j)x}-
	\sum_{i=1}^{n}\sum_{j=1}^{n}d_je^{2d_i}e^{(d_i+d_j)x}],\\
	&=\frac{1}
	{(\sum_{i=1}^{n}e^{d_ix})^2}
	[\sum_{i=1}^{n}\sum_{j=1}^{n}(d_i-d_j)e^{2d_i}e^{(d_i+d_j)x}].
	\end{aligned}
	$$
	Denote
	$$h(x):=
	\sum_{i=1}^{n}\sum_{j=1}^{n}(d_i-d_j)e^{2d_i}e^{(d_i+d_j)x}.$$
	It follows that $f^{'}(x)=\frac{h(x)}{(\sum_{i=1}^{n}e^{d_ix})^2}$.  Let
	$t_{ij}=e^{2d_i}e^{(d_i+d_j)x}$; then we have
	\begin{equation}\label{eq:3.1}
	h(x)=\sum_{i=1}^{n}\sum_{j=1}^{n}(d_i-d_j)t_{ij}.
	\end{equation}
	Exchanging the role of $i$ and $j$, we obtain
	\begin{equation}\label{eq:3.2}
	h(x)=\sum_{i=1}^{n}\sum_{j=1}^{n}(d_j-d_i)t_{ji}.
	\end{equation}
	Based on ($\ref{eq:3.1}$) and ($\ref{eq:3.2}$), we deduce that
	$$
	\begin{aligned}
	2h(x)
	= \sum_{i=1}^{n}\sum_{j=1}^{n}(d_i-d_j)(t_{ij}-t_{ji})
	=
	\sum_{i=1}^{n}\sum_{j=1}^{n}(d_i-d_j)(e^{2d_i}-e^{2d_j})e^{(d_i+d_j)x}.
	\end{aligned}
	$$
	It follows that $h(x)\geq0$ and hence $f^{'}(x)\geq0,\ \forall x\in(0,\infty)$. Therefore, $f(x)$ is a monotonic increasing function. The proof is completed.	
\end{proof}	
\end{appendices}

\bibliography{ref}% common bib file

%\input sn-sample-bib.tex%

\end{document}